\documentclass[draft,a4paper,12pt]{article}
\usepackage{amsfonts}
\usepackage{latexsym}
\usepackage{theorem}

\newcommand{\linespacing}{1.25}
\let\Section=\section
\newtheorem{theorem}{Theorem}[section]
\newtheorem{lemma}[theorem]{Lemma}
\newtheorem{proposition}[theorem]{Proposition}
\newtheorem{corollary}[theorem]{Corollary}

{\theorembodyfont{\normalfont\rmfamily}
\newtheorem{definition}[theorem]{Definition}
\newtheorem{remark}[theorem]{Remark}
}
\makeatletter
\def\ack{\vspace{.5\baselineskip}\noindent{\theorem@headerfont
Acknowledgement}\ \ }
\newenvironment{introthm}[1]%
 {\renewcommand{\theththm}{\ref{#1}}\begin{ththm}}
 {\end{ththm}}
\newtheorem{ththm}{Theorem}
\newenvironment{proof}[1][]%
{\def\proof@temp{#1}\par\noindent
\textsc{Proof}\ifx\proof@temp\@empty\else\ ({\proof@temp})\fi\hspace{1em}}
{\hphantom{xxx}\hfill~ {$\Box$}\par\vspace{.4\baselineskip}}
\def\operatorname#1{\mathop{\operator@font #1}\nolimits}%
\makeatother

\renewcommand{\L}{\mathcal{L}}
\newcommand{\C}{\mathbb{C}}
\newcommand{\R}{\mathbb{R}}
\newcommand{\Rnn}{\R^{2n}}

\newcommand{\Id}{\operatorname{Id}}

\newcommand{\supp}{\operatorname{supp}}

\newcommand{\bbnu}{[\![\nu]\!]}

\newcommand{\Cech}{\v{C}ech}

\def\ftnote#1{\def\footnotemark{}\footnote{#1}\setcounter{footnote}{0}}

\newdimen{\addresswidth}
\def\address#1#2{\parbox[t]{\addresswidth}{\centering
{{\small\ttfamily #1}\\[10pt]
\small{#2}}}}
\title{Traces for star products\\
on symplectic manifolds\ftnote{This research was partially supported by an
Action de Recherche Concert\'ee de la Communaut\'e fran\c{c}aise de
Belgique.}\\\ }
\author{%
Simone Gutt\\[10pt]
\address{sgutt@ulb.ac.be}
{Universit\'e Libre de Bruxelles\\
Campus Plaine, CP 218\\
BE -- 1050~Brussels\\
Belgium\\[5pt]and\\[5pt]
Universit\'e de Metz\\
Ile du Saulcy\\
57045 Metz Cedex 01\\
France}
\\[80pt]\\
John Rawnsley\\[10pt]
\address{J.Rawnsley@warwick.ac.uk}
{Mathematics Institute\\
University of Warwick\\
Coventry CV4 7AL\\
United Kingdom\\}
\\\ \\\ \\\
}

\date{May 2001}
\renewcommand{\baselinestretch}{\linespacing}

\begin{document}

\renewcommand{\baselinestretch}{1}

\maketitle

\thispagestyle{empty}

\begin{abstract}
We give a direct elementary proof of the existence of traces for
arbitrary star products on a symplectic manifold. We follow the approach
we used in \cite{refs:GuttRaw}, solving first the local problem. A
normalisation introduced by Karabegov \cite{refs:Karabegov} makes the
local solutions unique and allows them to be pieced together to solve
the global problem.
\end{abstract}

\newpage

\renewcommand{\baselinestretch}{\linespacing}

\Section{Introduction}\label{sect:intro}

In a previous paper \cite{refs:GuttRaw}, we gave lowbrow proofs of some
properties of differential star products on symplectic manifolds (in
particular the classification of equivalence classes of such star
products) using \Cech\ cohomology methods and the existence of special
local derivations which we called $\nu$-Euler derivations in
\cite{refs:GuttRaw}. In this note, we present, in a similar spirit,
properties of existence and uniqueness of traces for such star products.
Our proof of the existence of a trace relies on a canonical way of
normalization of the trace introduced by Karabegov
\cite{refs:Karabegov}, using local $\nu$-Euler derivations.

Let $*$ be a star product (which we always assume here to be defined by
bidifferential operators) on a symplectic manifold $(M,\omega)$. In the
algebra of smooth functions on $M$, consider the ideal $C^\infty_0(M)$
of compactly supported functions. A \textit{trace} is a $\C\bbnu$-linear
map $\tau \colon C^\infty_0 (M)\bbnu \to \C[\nu^{-1},\nu]\!]$ satisfying
\[
\tau(u*v) = \tau(v*u).
\]

The question of existence and uniqueness of such traces has been solved
by the following result.

\begin{introthm}{traces:exist}
{\rm(Fedosov \cite{refs:Fedosov1,refs:Fedosov2}; Nest--Tsygan
\cite{refs:NestTsygan})} Any star product on a symplectic manifold
$(M,\omega)$  has a trace which is unique up to multiplication by an
element of $\C[\nu^{-1},\nu]\!]$. Every trace is given by a smooth
density $\rho \in C^\infty(M)[\nu^{-1},\nu]\!]$:
\[
\tau(u) = \int_M u \rho \frac{\omega^n}{\nu^nn!}.
\]
\end{introthm}

We shall give here an elementary proof of this theorem. The methods
use intrinsically that we have a symplectic manifold. For Poisson manifolds
Felder and Shoikhet have shown in \cite{refs:FS} that the Kontsevich
star product also has a trace.

We are grateful to Martin Bordemann, Alexander Karabegov and Stefan Waldmann
for helpful comments.

\Section{Traces}\label{sect:traces}

Let $(M,\omega)$ be a connected symplectic manifold, $N$ the algebra of
smooth functions and $N_c$ the ideal in $N$ of compactly supported
functions. Obviously, $N_c\bbnu$ is an ideal in $N\bbnu$ and any
differential star product or equivalence on $N\bbnu$ is determined on
$N_c\bbnu$.

\begin{definition}
Let $*$ be a star product on $(M,\omega)$ then a \textit{trace} is
a $\C\bbnu$-linear map $\tau \colon N_c\bbnu \to \C[\nu^{-1},\nu]\!]$
satisfying
\[
\tau(u*v) = \tau(v*u).
\]
\end{definition}

\begin{remark}
Since any $\C[\nu^{-1},\nu]\!]$-multiple of a trace is a trace, it is
not necessary to work with Laurent series, but we do so for two reasons.
Firstly, in the formula for the trace of a pseudo-differential operator
a factor of $\nu^{-n}$ occurs, and secondly, the presence of such a
factor simplifies equation (\ref{traces:MoyalD}) below.
\end{remark}

If $\tau$ is a non-trivial trace, we take $u$ in $N_c$ and can then
expand
\[
\tau(u) = \nu^r \sum_{s\ge0} \nu^s \tau_s(u)
\]
where each $\tau_s \colon N_c \to \R$ is a linear map and we assume
$\tau_0 \ne 0$. The condition to be a trace takes the form
\begin{equation}\label{traces:trk}
\tau_k(\{u,v\}) + \tau_{k-1}(C_2^-(u,v)) + \dots +
\tau_0(C_{k+1}^-(u,v)) = 0
\end{equation}
for $k= 0, 1,2,\dots$,
where $C_r^-$ denotes the antisymmetric part of $C_r$.

\begin{remark}
In \cite{refs:CFS} the notion of a \textit{closed star product} was
introduced and related to cyclic cohomology of the algebra of functions.
The existence of a closed star product was proved in \cite{refs:OMY};
see also \cite{refs:BNW,refs:Pflaum}. A star product is closed if
\begin{equation}\label{traces:closed}
\int_M C_r^-(u,v) \omega^n = 0,\qquad \forall u,v \in N_c
\end{equation}
for all $1\le r\le n$ where $2n$ is the dimension of $M$. If
(\ref{traces:closed}) holds for all $r \ge 1$ then $u*v$ is said to be
\textit{strongly closed}. Thus to be strongly closed is the same as
requiring that
\[
\tau(u) = \int_M u \omega^n
\]
be a trace on $N_c\bbnu$.
\end{remark}

For $k=0$, equation (\ref{traces:trk}) reduces to the condition
$\tau_0(\{u,v\})=0$ for all $u$, $v$. In \cite{refs:BRW} it is shown by
elementary means that this implies that $\tau_0$ is a multiple of the
integral $\int_M u \omega^n$ when $M$ is connected. We give a proof here
for completeness.

\begin{lemma} {\rm(Gelfand \& Shilov \cite{refs:GS})}
If $u$ is a compactly supported smooth function on $\R^N$ with
$\int_{\R^N} u \, d^Nx = 0$ then $u$ is a sum of derivatives
of compactly supported smooth functions.
\end{lemma}

\begin{proof}
For $N=1$ we simply set $v(x) = \int_{-\infty}^x u(t) dt$ and observe that
$v$ obviously vanishes when $x$ is below the support of $u$, and is zero
again when $x$ is large as a consequence of $\int_{-\infty}^\infty u(t) dt
=0$. Thus $v$ is compactly supported and $u(x) = \frac{dv}{dx}$. Moreover,
if $u$ depends smoothly on parameters, so will $v$ and if $u$ is compactly
supported in the parameters, so is $v$.

Now proceed by induction on $N$. $w(x_1,\ldots,x_{N-1}) =
\int_{-\infty}^\infty u(x_1,\ldots,x_{N-1}, t)\, dt$ clearly is
compactly supported in $\R^{N-1}$ and has vanishing integral, so by the
inductive assumption $w(x_1,\ldots,x_{N-1}) =
\displaystyle{\sum_{i=1}^{N-1}\frac{\partial
w_i(x_1,\ldots,x_{N-1})}{\partial x_i}}$ for some compactly supported
functions $w_i$ on $\R^{N-1}$. Take a compactly supported bump function
$r(t)$ on $\R$ with $\int_{\R} r(t) \, dt=1$ and consider
$u(x_1,\ldots,x_N) - w(x_1,\ldots,x_{N-1})r(x_N)$ which is compactly
supported in all its variables. Integrating in $x_N$ we see that the
integral over $\R$ vanishes, and so $u(x_1,\ldots,x_N) -
w(x_1,\ldots,x_{N-1})r(x_N) = \displaystyle{\frac{\partial v}{\partial
x_N}}(x_1,\ldots,x_N)$. Thus $ u = \displaystyle{\frac{\partial
v}{\partial x_N}}(x_1,\ldots,x_N) +
\displaystyle{\sum_{i=1}^{N-1}\frac{\partial w_i(x_1,\ldots,x_{N-1})
r(x_N)}{\partial x_i}}$ which completes the inductive step.
\end{proof}

\begin{lemma} {\rm(Bordemann, R\"omer, Waldmann \cite{refs:BRW})}
Let $(M,\omega)$ be a connected symplectic manifold. If $\sigma \colon
N_c \to \R$ is a linear map with $\sigma(\{u,v\}) = 0$  for all $u,v \in
N_c$ then $\sigma(u) = c \int_{M} u \omega^n$ for some constant $c$.
\end{lemma}

\begin{proof}
Fix $u \in N_c$ and cover $M$ by Darboux charts $U_\alpha$ such that
only finitely many $U_\alpha$ intersect the support of $u$. Take
a partition of unity $\varphi_\alpha$ subordinate to $U_\alpha$ then
only a finite number of $u\varphi_\alpha$ are non-zero. Thus
$\sigma(u) = \sum_\alpha \sigma(u\varphi_\alpha)$.

$u\varphi_\alpha - \displaystyle{\frac{\int_M u\varphi_\alpha \omega^n}{\int_M
\varphi_\alpha \omega^n}} \varphi_\alpha$ has vanishing integral on
$U_\alpha$ which can be viewed as an open set in some $\R^{2n}$. Thus by
the previous Lemma there are functions $v_i$, $w_i$ with compact support
such that $u\varphi_\alpha - \displaystyle{\frac{\int_M u\varphi_\alpha
\omega^n}{\int_M \varphi_\alpha \omega^n}} \varphi_\alpha = \sum_i
\displaystyle{\frac{\partial v_i}{\partial p_i}} + 
\displaystyle{\frac{\partial w_i}{\partial q_i}}$
for the Darboux coordinates $p_i$ and $q_i$. But $\displaystyle{\frac{\partial
v_i}{\partial p_i}} = \{v_i, q_i\}$, and if we choose a function $s_i$ of
compact support which is identically 1 on the support of $v_i$ then
$\displaystyle{\frac{\partial v_i}{\partial p_i}} = \{v_i, s_iq_i\}$ so we see that $
\sigma\left(\displaystyle{\frac{\partial v_i}{\partial p_i}}\right) =0$ and similarly $
\sigma\left(\displaystyle{\frac{\partial w_i}{\partial q_i}}\right) = 0$. Thus
$\sigma(u\varphi_\alpha) = \displaystyle{\frac{\int_M u\varphi_\alpha \omega^n}{\int_M
\varphi_\alpha \omega^n}} \sigma(\varphi_\alpha)$. Hence
\begin{eqnarray*}
\sigma(u) &=& \int_M u\sum_\alpha\displaystyle{\frac{ \varphi_\alpha }{\int_M
\varphi_\alpha \omega^n}} \sigma(\varphi_\alpha)\omega^n\\
&=& \int_M u \rho \omega^n.
\end{eqnarray*}
Since $\sigma(\{u,v\}) = 0$, $\int_M \{u,v\} \rho \omega^n = 0$
and hence $\int_M \{\rho,u\} v \omega^n = 0$ for all $v$. Thus
$\{\rho,u\} = 0$ for all $u$ and hence $\rho = c$, a constant.
\end{proof}

Thus any trace has
the form
\[
\tau(u) = a\nu^r\left(\int_M u \omega^n + \sum_{k\ge1} \nu^k \tau_k(u)\right)
\]
where $a \ne 0$. We can divide by $a$ and multiply by $\nu^{-n-r}$
to bring $\tau$ into the form
\[
\int_M u \frac{\omega^n}{\nu^n} + \nu^{-n}\sum_{k\ge1} \nu^k \tau_k(u).
\]
Any trace in this form will be said to be \textit{standard} c.f.\
\cite{refs:Guillemin}.

As observed in \cite{refs:BRW}, this is enough to show that any two
traces $\tau$ and $\tau'$ for the same star product are proportional.
For, if  $\tau$ is standard, the leading term of $\tau'$ is a multiple
$c \nu^{r}$ of the integral, hence is equal to the leading term of
$\tau$ multiplied by $c \nu^{r+n}$ . But then $\tau' - c\nu^{r+n}\tau$
is a trace which vanishes to at least order $r+1$ in $\nu$. This
argument can be repeated indefinitely to show that $\tau' =
c\nu^{r+n}(1+ \sum_k \nu^k c_k) \tau$. Remark that in particular if
$\tau$ and $\tau'$ are standard then $\tau' = (1+ \sum_k \nu^k c_k)
\tau$. This proves

\begin{theorem} {\rm(Nest \& Tsygan \cite{refs:NestTsygan})}\label{thm:prop} 
On a connected symplectic manifold $(M,\omega)$ any two traces are
proportional by an element of $\C[\nu^{-1},\nu]\!]$.
\end{theorem}

\Section{The local case}\label{sect:local}

In the case of $\Rnn$ with its standard constant $2$-form $\Omega$,
then
\[
\tau_M(u) = \int_{\Rnn} u\, \textstyle{\frac{\Omega^{n}}{\nu^n n!}}
\]
is a trace on compactly supported functions for the Moyal star product
$*_M$. The Moyal star product and this trace have an important
homogeneity property \cite{refs:Karabegov}. If we take a conformal
vector field $\xi$ on $\Rnn$, so $ \L_\xi\Omega = \Omega$ then $D_M =
\xi + \nu \frac{\partial}{\partial\nu}$ is a derivation of $*_M$ and
$\tau_M$ satisfies
\begin{equation}\label{traces:MoyalD}
\tau_M(D_Mu) = \nu \frac{\partial}{\partial\nu} \tau_M(u).
\end{equation}

If $*$ is any star product defined on an open ball $U$ in $\Rnn$ then it
is equivalent to the restriction of the Moyal star product by a map $T =
\Id + \sum_{k\ge1} \nu^k T_k$ with
\[
T(u*v) = T(u)*_M T(v)
\]
and then we see that
\[
\tau (u) = \int_U T(u) \frac{\Omega^n}{\nu^{n}}
\]
is a trace for $*$. Each $T_k$ is a differential operator, so it has a
formal adjoint $T'_k$ so that, if we put $T' = \Id + \sum_{k\ge1} \nu^k
T'_k$, then
\[
\tau (u) =\int_U u T'(1) \frac{\Omega^n}{\nu^{n}}
\]
for $u \in C^{\infty}_c(U)$. If we put $\rho = T'(1) \in
C^{\infty}(U)[\nu^{-1},\nu]\!]$ then
\[
\tau (u) = \int_U u \rho \frac{\Omega^n}{\nu^{n}}.
\]
If $\tau$ is standard then $\rho = 1 + \sum_{k\ge1} \nu^k \rho_k$.

Further $D = T^{-1} \circ D_M \circ T$ will be a derivation of $*$
and satisfies the transform of equation (\ref{traces:MoyalD}):
\[
\tau(Du) = \nu \frac{\partial}{\partial\nu} \tau(u).
\]
$D$ has the form $\xi + \nu \frac{\partial}{\partial\nu} + D'$. Local
derivations of this form we give a special name

\begin{definition}
Let $(M,\omega)$ be a symplectic manifold. Say that a derivation $D$ on
an open set $U$ of $N\bbnu,*$ is \textbf{$\nu$-Euler} if it has the form
\begin{equation}\label{intrinsic:dereqn}
D = \nu \frac{\partial}{\partial\nu} + X + D'
\end{equation}
where $X$ is conformally symplectic ($\L_X\omega=\omega$) and
$D' = \sum_{r\ge1} \nu^r D'_r$ with the $D'_r$ differential
operators on $U$.
\end{definition}

Note that conformally symplectic vector fields only exist locally in
general, so we also cannot ask for global $\nu$-Euler derivations. Two
local $\nu$-Euler derivations defined on the same open set will differ
by a $\nu$-linear derivation, and so the difference is
$\nu^{-1}\times{}$inner. This means that $\tau \circ D$ will be
independent of $D$ and thus is globally defined as an $\R$-linear
functional, even if $D$ is not.

\begin{definition}
A standard trace $\tau$ is \textit{normalised} if it satisfies the
analogue of the Moyal homogeneity condition:
\[
\tau(Du) = \nu\frac{\partial}{\partial\nu} \tau(u)
\]
on any open set where there are $\nu$-Euler derivations and for any
such local $\nu$-Euler derivation $D$.
\end{definition}

This condition was introduced by Karabegov \cite{refs:Karabegov}.

It is clear that the pull-back of the Moyal trace by an equivalence
with the Moyal star product is normalised, so a normalised trace always
exists for any star product on an open ball in $\Rnn$.

\begin{proposition} {\rm(Karabegov \cite{refs:Karabegov})}
If $\tau$ and $\tau'$ are normalised traces for the same star
product on an open ball $U$ in $\Rnn$ then $\tau=\tau'$.
\end{proposition}

\begin{proof} Neither trace can be zero, and are proportional so
$\tau' = (1+c \nu^r + \dots)\tau$. If $\tau'\ne\tau$ then there is a first
$r>0$
where $c\ne0$. Then we substitute in the normalisation condition to give
\begin{eqnarray*}
(1+c \nu^r + \dots)\tau(Du) &=& \tau'(Du)\\
&=&  \nu \frac{\partial\ }{\partial\nu} \tau'(u)\\
&=&  \nu \frac{\partial\ }{\partial\nu}((1+c \nu^r + \dots)\tau(u))\\
&=& (rc\nu^r+\dots)\tau(u) + (1+c \nu^r + \dots)\tau(Du)
\end{eqnarray*}
which implies that $c=0$. This contradiction shows that $\tau'=\tau$.
\end{proof}

\Section{The global case}\label{sect:global}

Let $(M,\omega)$ be a connected symplectic manifold. Then we can cover
$M$ by Darboux charts $U$ which are diffeomorphic to open balls in
$\Rnn$ and such that all non-empty intersections are also diffeomorphic
to open balls. Let $*$ be a star product on $M$ and then the restriction
of $*$ to $U$ has a normalised trace with density $\rho_U \in
C^\infty(U)\bbnu$. If we have two open sets $U$, $V$ which overlap, then
on the intersection both $\rho_U$ and $\rho_V$ will determine normalised
traces on an open ball in $\Rnn$. Since there is only one such trace,
$\rho_U = \rho_V$ on $U \cap V$. It follows that there is a globally
defined function $\rho$ on $M$ such that $\rho_U = \rho|_U$. Set
\[
\tau(u) = \int_M u \rho \frac{\omega^n}{\nu^n}
\]
for $u \in C^\infty_c(M)$.

Given $u$, $v$ in $C^\infty_c(M)$ we can find a finite partition of
unity $\varphi_i$ on $\supp u \cup \supp v$ with supports in the open
sets above. Then $u=\sum_i \varphi_i u$ and $v=\sum_j \varphi_j v$ so
$u*v = \sum_{i,j} (\varphi_i u)*(\varphi_j v)$ so $\tau(u*v) =
\sum_{i,j} \tau((\varphi_i u)*(\varphi_j v))= \sum_{i,j}\tau((\varphi_j
v)*(\varphi_i u)) = \tau(v*u)$ so $\tau$ is a trace. A similar partition
of unity argument shows that $\tau$ is normalised.

Combining this with Theorem \ref{thm:prop} and the fact that the
normalised trace we just constructed has a smooth density we obtain:

\begin{theorem} [Fedosov, Nest--Tsygan]\label{traces:exist}
On a connected symplectic manifold $(M,\omega)$ any differential star
product has a unique normalised trace. Any trace is multiple of this
and is given by a smooth density.
\end{theorem}

\begin{remark}
In \cite{refs:NestTsygan} a proof that any trace has a smooth density is
given using cyclic cohomology.
\end{remark}

\begin{corollary}
Any trace is invariant under all $\C\bbnu$-linear
automorphisms of the star product.
\end{corollary}

\begin{proof} A smooth trace is a multiple (in $\C[\nu^{-1},\nu]\!]$) of
a normalised trace. The transform of a $\nu$-Euler derivation by a
$\C\bbnu$-linear automorphism is again a $\nu$-Euler derivation, thus
the transform of a normalised trace by a $\C\bbnu$-linear automorphism
is again a normalised trace, and so is equal to the original normalised
trace.
\end{proof}

\end{document}